\documentclass[12pt]{amsart}
\usepackage{latexsym,cite,amsthm,amsmath,amssymb}

\usepackage{amsmath}
\usepackage{amsthm}
\date{}
\def\bc{\begin{center}}
\def\ec{\end{center}}
\def\phi{\varphi}
\def\epsilon{\varepsilon}

\def\l{\lambda}

\def\a{\alpha}

\def\0{\bar 0}
\def\1{\bar 1}

\def\F{{F\langle x,y,z\rangle}}
\def\SJ{{SJ[x,y,z]}}

\def\ctd{\hfill$\Box$}

\def\bes{\begin{eqnarray*}}
\def\ees{\end{eqnarray*}}
\def\bee{\begin{eqnarray}}
\def\eee{\end{eqnarray}}

\def\le{{\langle}}
\def\re{{\rangle}}
\def\Proof{{\it Proof. }}

\newtheorem{Th}{Theorem}
\newtheorem{Lem}{Lemma}

\newtheorem{Cor}{Corollary}

\begin{document}
\title{On commuting $U$-operators in Jordan algebras}
\vspace {5 mm}

\author {Ivan Shestakov}
\address{Instituto de Matem\'atica e Estat\'\i stica,
         Universidade de S\~ao Paulo,
                  S\~ao Paulo, Brazil,
                  and
                  Sobolev Institute of Mathematics, Novosibirsk, Russia}
         \email{shestak@ime.usp.br}
  \vspace {5 mm}

\maketitle
\begin{abstract}
Recently J.\,A.\,Anquela, T.\,Cort\'es, and H.\,Petersson \cite{ACP} proved that for elements $x, y$ in a non-degenerate  Jordan algebra $J$, the relation $x \circ y = 0$  implies that the $U$-operators of $x$  and $y$ commute: $U_xU_y = U_yU_x$. We show that the result may be not true without the assumption on non-degeneracity of $J$. We give also a more simple proof of the mentioned result in the case of linear Jordan algebras, that is, when $char\, F\neq 2$.
\end{abstract}

\medskip
\begin{flushright}
{Dedicated to Professor Amin Kaidi \\on the occasion of his 65-th annyversary}
\end{flushright}

\section{An Introduction}
\hspace{\parindent}

In a recent paper \cite{ACP} J.\,A.\,Anquela, T.\,Cort\'es, and H.\,Petersson have studied the following question for Jordan algebras:

(1) does the relation $x\circ y = 0$ imply that the quadratic operators $U_x$ and $U_y$ commute?

They proved that the answer is positive for non-degenerate Jordan algebras, and left open the question in the general case, not assuming nondegeneracy.

\smallskip

We show that the answer to question (1) is negative in general case. We give also a more simple proof of the result for linear non-degenerate Jordan algebras, that is, over a field $F$ of characteristic $\neq 2$.

Unless otherwise stated, we will deal with associative and Jordan algebras over a field of arbitrary characteristic.

\section{A counter-example}
\hspace{\parindent}

Let us recall some facts on Jordan algebras. We use as general  references the books  \cite{Jac, ZSSS,  McC}, and the paper \cite{McC1}.

\smallskip
Consider the free special Jordan algebra $SJ[x,y,z]$ and the free associative algebra $F\langle x,y,z\rangle$ over a field $F$. Let $*$ be the involution of $\F$ identical on the set $\{x,y,z\}$. Denote $\{u\}=u+u^*$ for $u\in \F$, then $\{u\}\in\SJ$ \cite{Jac, ZSSS} (see also \cite{McC1} for the case of characteristic 2). Below $ab$ will denote the associative product  in $\F$, so that  $a\circ b=ab+ba$ and  $aU_b=bab$ are the corresponding  linear and quadratic operations  in $\SJ$. 

\smallskip
For an ideal $I$ of $SJ[x,y,z]$, let $\hat I$ denote the ideal of  $F\langle x,y,z\rangle$ generated by $I$. By Cohn's Lemma \cite[lemma 1.1]{Jac} (see also \cite[Corollary to Cohn's Criterion]{McC1}), the quotient algebra  $J=\SJ/I$ is special if and only if $I=\hat{I}\cap \SJ$.

\begin{Lem}\label{z[Ux,Uy]} The following equality holds in $\SJ\subseteq \F$:
$$
z[U_x,U_y]=\{(x\circ y)zxy\}-zU_{x\circ y}.
$$
\end{Lem}
\Proof We have in $\F$
\bes
 z[U_x,U_y]&=&yxzxy-xyzyx=(y\circ x)zxy-xyzxy-xyzyx=\\
&=&(y\circ x)zxy-xyz(x\circ y)=\{(x\circ y)zxy\}-(x\circ y)z(x\circ y).
\ees 
\ctd
 \begin{Th} Let $I$ denote the ideal of $\SJ$ generated by $x\circ y=xy+yx$ and $J=\SJ/I$. Then for the images $\bar x,\bar y$ of the elements $x,y$ in $J$ we have $\bar x \circ \bar y=0$ but $[U_{\bar x},U_{\bar y}]\neq 0$.
\end{Th}  
 \Proof It suffices to show that $k=z[U_x,U_y]\notin I$. By lemma 1, $k=\{(x\circ y)zxy\}\ (mod\ I)$. Now, the arguments from the proof of \cite[theorem 1.2]{Jac}, show that $k\notin I$ when $F$ is a field of characteristic not 2 (see also \cite[exercise 1, page 12]{Jac}). 

\smallskip

The result is also true  in characteristic 2 for quadratic Jordan algebras. 
In this case, one needs certain modifications concerning the generation of ideals in quadratic case. The author is grateful to T.\, Cort\'es and J.\,A.\, Anquela who 
corrected the first ``naive'' author's proof and suggested the proper modifications which we give below.

\smallskip
We have to prove  that  $\{(x \circ y)zxy\} \not\in I$. By \cite[(1.9)]{NMcC}, the ideal $I$ is the outer hull of $F(x\circ y)+U_{x\circ y} \widehat{SJ[x,y,z]}$, where $\widehat J$ denotes the unital hull of $J$. Assume that there exists a Jordan polynomial $f(x, y, z, t) \in SJ[x, y, z, t]$ with all of its Jordan monomials containing the variable $t$, such that $\{(x \circ y)zxy\} = f(x, y, z, x \circ y)$. By degree considerations, $f = g + h$, where $g, h \in SJ[x, y, z, t]$, $g$ is multilinear, and $h(x, y, z, t)$ is a linear combination of $U_t z$ and $z\circ t^2$. On the other hand, arguing as in \cite[Theorem 1.2]{Jac}, $g \in SJ[x, y, z, t] \subseteq H(F\langle x, y, z, t\rangle, *)$, and because of degree considerations and the fact that $z$ occupies inside position in the associative monomials of $\{(x\circ y)zxy\}$, $g$ is a linear combination of
$$
\{xzyt\}, \{xzty\}, \{tzxy\}, \{tzyx\}, \{yztx\}, \{yzxt\},
$$
and $h$ is a scalar multiple of $U_t z$. Hence $f$ has the form
\bes
f(x, y, z, t) &=& \a_1\{xzyt\} + \a_2\{xzty\} + \a_3\{tzxy\}\\
&+& \a_4\{tzyx\} + \a_5\{yztx\} + \a_6\{yzxt\}\\
&+& \a_7 tzt,
\ees
and therefore
\bes
\{(x \circ y)zxy\} &=& \a_1\{xzy(x\circ y)\} + \a_2\{xz(x\circ y)y\} + \a_3\{(x\circ y)zxy\}\\
&+& \a_4\{(x\circ y)zyx\} + \a_5\{yz(x\circ y)x\} + \a_6\{yzx(x\circ y)\}\\
&+& \a_7 (x\circ y)z(x\circ y),
\ees
Comparing coefficients as in \cite[Theorem 1.2]{Jac}, we get 
\bes
\a_1 = \a_2 = \a_5 = \a_6 = 0,\\
\a_3 = \l + 1, \a_4 = \l, \a_7 = - 2\l,
\ees
for some $\l\in F$. Going back to $f$, we get
$$
f = (\l + 1)\{tzxy\} + \l\{tzyx\} - 2 \l tzt = \{tzxy\} + \l\{tz(x\circ y)\} - 2\l U_t z,
$$
so that $\{tzxy\} \in SJ[x, y, z, t]$, which is a contradiction.

 In fact, the standard arguments with the Grassmann algebra do not work in characteristic 2, to prove that $\{tzxy\} \notin SJ[x, y, z, t]$, but one can check directly (or with aid of computer) that the space of symmetric multilinear elements in $F\le x,y,z,t\re$ has dimension 12 while the similar space of Jordan elements has dimension 11.

\ctd 
 
\section{The non-degenerate case}
\hspace{\parindent}

Here we will give another proof of the main result from \cite{ACP} that the answer to question (1) is positive for nondegenerate algebras, in the case of  linear  Jordan algebras (over a field $F$ of characteristic $\neq 2$).

Let J be a linear Jordan algebra, $a\in J,\ R_a : x \mapsto xa$ be the operator of right multiplication on $a$, and $U_a = 2R^2_a- R_{a^2}$.

As in \cite{ACP}, due to the McCrimmon-Zelmanov theorem \cite{MZ}, it suffices to consider Albert algebras. We will need only the fact that an Albert algebra $A$ is {\em cubic}, that is, for every $a\in A$, holds the identity
\bes
a^3=t(a)a^2-s(a)a+n(a),
\ees
where $t(a), s(a), n(a)$ are linear, quadratic, and cubic forms  on $A$, correspondingly \cite{Jac}. Linearizing the above identity on $a$, we get the identity
\bes
2((ab)c+(ac)b+(bc)a)&=&2(t(a)bc+t(b)ac+t(c)ab)\\
&&-s(a,b)c-s(a,c)b-s(b,c)a+n(a,b,c),
\ees
where $s(a,b)=s(a+b)-s(a)-s(b)$ and $n(a,b,c)=n(a+b+c)-n(a+b)-n(a+c)-n(b+c)+n(a)+n(b)+n(c)$ are bilinear and trilinear forms. In particular, we have
  \bee\label{aab}
a^2b+2(ab)a=t(b)a^2+2t(a)ab-s(a,b)a-s(a)b+\tfrac12n(a,a,b).
\eee

\begin{Lem}\label{lem[Ua,Ub]}
Let $a,b\in J$ with $ab=0$. Then $[U_a,U_b]=[R_{a^2},R_{b^2}]$.
\end{Lem}  
\Proof
Linearizing the Jordan identity $[R_x,R_{x^2}]=0$, one obtains
$$
[R_{a^2},R_b]=-2[R_{ab},R_a]=0,
$$
and similarly $[R_a,R_{b^2}]=0$. Therefore,
\bes
[U_a,U_b]=[2R_a^2-R_{a^2},2R_b^2-R_{b^2}]=4[R_a^2,R_b^2]+[R_{a^2},R_{b^2}].
\ees
Furthermore, $[R_a^2,R_b^2]=[R_a,R_b^2R_a+R_aR_b^2]$. By the operator Jordan identity \cite[(1.$O_2$)]{Jac},
\bes
R_b^2R_a+R_aR_b^2=-R_{(ba)b}+2R_{ab}R_b+R_{b^2}R_a=R_{b^2}R_a,
\ees
therefore $[R_a^2,R_b^2]=[R_a,R_{b^2}R_a]=[R_a,R_{b^2}]R_a=0$, which proves the lemma. \ctd
 
\begin{Th}
Let $J$ be a cubic Jordan algebra and $a,b\in J$ with $ab=0$. Then $[U_a,U_b]=0$.
\end{Th}
\Proof For any $c\in J$ we have by Lemma 2 and by the linearization of the Jordan identity $(x,y,x^2)=0$
\bes
c[U_a,U_b]=c[R_{a^2},R_{b^2}]=(a^2,c,b^2)=-2(a^2b,c,b).
\ees
By (\ref{aab}), we have
\bes
(a^2b,c,b)&=&t(b)(a^2,c,b)-s(a)(b,c,b)-s(a,b)(a,c,b)\\
&=&-2t(b)(ab,c,a)-s(a,b)(a,c,b)=-s(a,b)(a,c,b).
\ees
Substituting $c=a$, we get $(a^2b,a,b)=((a^2b)a)b=(a^2(ba))b=0$, which implies
$0=s(a,b)(a,a,b)=s(a,b)(a^2b)$. Therefore, $s(a,b)=0$ or $a^2b=0$. In both cases
this implies $c[U_a,U_b]=0$.
\ctd
\begin{Cor}
In an Albert algebra $A$, the equality $ab=0$ implies $[U_a,U_b]=0$.
\end{Cor}

In connection with the counter-example above, we would like to formulate an open question. Let $f,g\in SJ[x,y,z]$ such that $g\in \widehat{(f)}$ but $g\not\in (f)$, where $(f)$ and  $\widehat{(f)}$ are the ideals generated by $f$ in $SJ[x,y,z]$ and in $F\le x,y,z\re$, respectively. Then the quotient algebra $SJ[x,y,z]/(f)$ is not special, due to Cohn's Lemma. It follows from the results of \cite{Z} that the quotient algebra $\widehat{(f)}/(f)$ is degenerated. The question we want to ask is the following: 
\smallskip

{\em If $f=0$ in a nondegenerate Jordan algebra $J$, should also be $g=0$?}
\smallskip

Of course, there is a problem of writing $f$ and $g$ in an arbitrary Jordan algebra,  we know only what they are in $SJ[x,y,z]$, but in the free Jordan algebra $J[x,y,z]$ they have many pre-images (up to $s$-identities), and one may choose pre-images for which the question has a negative answer. For example, the answer is probably negative for $f=x\circ y$ and $g=z[U_x,U_y]+G(x,y,z)$, where $G(x,y,z)$ is the Glennie $s$-identity \cite{Jac}.  

So we modify our question in the following way: 

{\it In the situation as above, is it true that there exists $g'\in J[x,y,z]$ such that $g - g'$ is an $s$-identity and $f=0$ implies $g'=0$ in non-degenerate Jordan algebras?}

\section{Acknowledgements}
The author acknowledges the support by FAPESP, Proc.\,2014/09310-5 and CNPq, Proc.\,303916/ 2014-1. He is grateful to professor Holger Petersson for useful comments and suggestions, and to professors Jos\'e \'Angel Anquela and Teresa Cort\'es for  correction  the proof of Theorem 1 in the case of characteristic 2. He thanks all of them for pointing out some misprints.

\end{document}